\documentclass[11pt, a4paper]{amsart}
\usepackage{fontenc, amsmath ,mathrsfs,amsfonts, amsthm, geroENG, amssymb, MnSymbol, multicol, enumerate}
\usepackage{color}

\usepackage[english]{babel}
\usepackage{hyperref}
\usepackage{garamondlibre} 
\usepackage{mathrsfs}
\begin{document}
\date{}

\title{A Complex Borel-Bernstein Theorem}
\author{Gerardo González Robert}
\address{Facultad de Ciencias, Universidad Nacional Autónoma de México \\ Ciudad de México, México}
\email{matematicas.gero@gmail.com, gero@ciencias.unam.mx}

\begin{abstract}
Zero-one laws are a central topic in metric Diophantine approximation. A classical example of such laws is the Borel-Bernstein theorem. In this note, we prove a complex analogue of the Borel-Bernstein theorem for complex Hurwitz continued fractions. As a corollary, we obtain a complex version of Khinchin's theorem on Diophantine approximation.
\end{abstract}

\maketitle
\section{Introduction}
For over a century, number theory and probability have had a fruitful relationship. A particularly fecund environment for this alliance is the theory of regular continued fractions, which we shall denote by
\[
[a_0;a_1,a_2,\ldots]_{\RE}= a_0 + \frac{1}{a_1 + \cfrac{1}{a_2 + \cfrac{1}{\ddots}}}.
\]
The Borel-Bernstein Theorem is a classical result in this area, its proof can be found in \cite{bugeaud-livre} (Theorem 1.11.). 

\begin{teo01}[Borel-Bernstein]\label{TeoBB}
Let $\bfu=(u_n)_{n\geq 1}$ be a sequence of positive real numbers and define 
\[
E_{\RE}(\bfu) := \left\{ z=[0;a_1,a_2,\ldots]_{\RE} \in [0,1] : a_n \geq u_{n}\; \text{ for infinitely many } n\in\Na\right\}.
\]
If $\leb_1$ denotes the Lebesgue measure on $[0,1)$, then
\[
\leb_1(E_{\RE}(\bfu)) = 
\begin{cases}
0, \text{ if } \displaystyle\sum_{n\in\Na} u_n^{-1} < +\infty, \\
1, \text{ if } \displaystyle\sum_{n\in\Na} u_n^{-1} = +\infty.
\end{cases}
\]
\end{teo01}
In this note, we show a complex analogue of the Borel-Bernstein theorem using Hurwitz continued fractions (defined on Section 2). Hurwitz continued fractions allow us to represent every complex number as a finite or infinite continued fraction of the form
\[
[a_0;a_1,a_2,a_3,\ldots] := a_0 + \cfrac{1}{a_1 + \cfrac{1}{a_2 + \cfrac{1}{\ddots}}},
\]
where each $a_n$ is a Gaussian integer. In this context, a unit square $\mfF$ centered at the origin (defined on Section 2) will play the role the interval $[0,1)$ plays within the regular continued fraction theory.

\begin{teo01}[Main result]\label{TEOBBC}
Let $\bfu=(u_n)_{n\geq 1}$ be a sequence of positive real numbers and define 
\[
E(\bfu) = \left\{ z=[0;a_1,a_2,\ldots]\in\mfF: |a_n|\geq u_{n}\; \text{ for infinitely many } n\in\Na\right\}.
\]
If $\leb$ is the Lebesgue measure on $\mfF$, then
\[
\leb(E(\bfu)) = 
\begin{cases}
0, \text{ if } \displaystyle\sum_{n\in\Na} u_n^{-2} < +\infty, \\
1, \text{ if } \displaystyle\sum_{n\in\Na} u_n^{-2} = +\infty.
\end{cases}
\]
\end{teo01}

We obtain a complex version of Khinchin's theorem on Diophantine approximation (Theorem 1.10 in \cite{bugeaud-livre}) as a consequence of Theorem \ref{TEOBBC}. 
\begin{coro01}\label{COROKHIN}
Given $\psi:\RE_{\geq 1} \to\RE_{>0}$, define the set
\[
K(\psi)= \left\{ \xi\in \Cx: \left| \xi-\frac{p}{q}\right| \leq  \psi(|q|)\; \text{ for infinitely many } p,q\in\Za[i], q\neq 0\right\}.
\]
If $x\mapsto x^2\psi(x)$ is non-increasing, then
\[
\leb(K(\psi))=
\begin{cases}
0, \text{ if } \displaystyle\sum_{n\geq 1} n^3\psi(n)^2 < +\infty,\\
\text{full}, \text{ if } \displaystyle\sum_{n\geq 1} n^3\psi(n)^2 = +\infty.
\end{cases}
\]
\end{coro01}
A. Nogueira studied a Borel-Bernstein theorem for some multidimensional continued fractions algorithms, such as the Perron-Jacobi and Brun algorithm \cite{nogueira}. Some results resembling Corollary \ref{COROKHIN} are well known. In 1952, W. J. Leveque introduced a complex continued fraction expansion and he used it to show a Khinchin-type result (Theorem 5 in \cite{leveque}). In 1982, with the aid of disjoint spheres and Bianchi groups, D. Sullivan obtained a result stronger than Corollary \ref{COROKHIN} (see Theorem 5 in \cite{sullivan}).  {In 2006, V. Beresnevich, D. Dickinson, and S. Velani strengthened Sullivan's result by applying their work on ubiquitous systems. They proved that Corollary \ref{COROKHIN} still holds if we replace ``$x\mapsto x^2\psi(x)$ is non-increasing'' by ``$\psi$ is non-increasing'' (Theorem 7 in \cite{bdv}). Moreover, in the recent paper \cite{HeXi}, Y. He and Y. Xiong used Hurwitz continued fractions to compute the Hausdorff dimension of $K(\psi)$ for non-increasing functions $\psi$ such that $\psi(x)$ tends to $0$ when $x$ tends to $+\infty$. We can combine Theorem 8 in \cite{bdv} and Theorem 1.2 in \cite{HeXi} to obtain the next result:
\begin{teo01}\label{Teo-BDV-HX}
Let $\psi:\RE_{>0}\to\RE_{>0}$ be a non-increasing function such that $\psi(x)\to 0$ when $x\to +\infty$. Define 
\[
\lambda:=\lambda(\psi) = \liminf_{x\to +\infty} \frac{-\log\psi(x)}{\log(x)}
\]
and let $\HD(\psi)$ be the Hausdorff dimension of $K(\psi)$. If $\lambda>2$, then $\HD(\psi)=\tfrac{4}{\lambda}$ and $\HD(\psi)=2$ whenever $\lambda\leq 2$. Furthermore, if $\clH^{\HD(\psi)}$ denotes the $\HD(\psi)$-Hausdorff measure on $\Cx$, we have
\[
\clH^{\HD(\psi)}(K(\psi))=
\begin{cases}
0, \text{ if } \sum_{n\geq 1} \psi(n)^{\HD({\psi}) } n^3<+\infty,\\
\text{full, if } \sum_{n\geq 1} \psi(n)^{\HD({\psi}) } n^3=+\infty.
\end{cases}
\]
\end{teo01}
When we choose $\veps>0$ and $\psi(x)=x^{-(1+\veps)}$, Theorem \ref{Teo-BDV-HX} is precisely the complex Jarník-Besicovitch Theorem (Theorem 6.1 in \cite{dokr}). A nice discussion of Khinchin-type results and a simple proof of the complex Jarník-Besicovitch Theorem can be found in \cite{dokr}.} A related multidimensional theorem was proven by M. Hussain in \cite{hussain17}. None of these papers use on any Borel-Bernstein-type result. However, they do rely on the convergence part of the Borel-Cantelli Lemma, whose proof can be found in almost every measure theory book.

\begin{lem01}[Borel-Cantelli lemma, convergence part]
Let $(X,\scA,\mu)$ be a measure space and $(A_n)_{n\geq 1}$ a sequence in $\scA$. If $\sum_n \mu(A_n)<+\infty$, then $\mu(\limsup_n A_n)=0$.
\end{lem01}

The organization of this note is as follows: in the next section, we define Hurwitz continued fractions and prove some useful lemmas; afterwards, we prove Theorem \ref{TEOBBC} and Corollary \ref{COROKHIN}.

\textbf{Notation.} The natural numbers, $\Na$, are the positive integers and $\Na_0:=\Na\cup\{0\}$. If $\saxu$ and $\sayu$ are two sequences of non-negative real numbers, $x_n\ll y_n$ means that there is some constant $\kappa>0$ such that $x_n\leq \kappa y_n$ holds for sufficiently large $n\in\Na$. We write $x_n\asymp y_n$ if $x_n\ll y_n$ and $y_n\ll x_n$ are true. Whenever the constant implied by $\ll$ depends on some object $\alpha$, we shall write $\ll_{\alpha}$. For any $A\subseteq \Cx$, $A^*:=A\setminus\{0\}$, $\Cl(A)$ is the closure of $A$, $\inte A$ is the interior of $A$, $zA:=\{za:a\in A\}$ for a given $z\in\Cx$, and, if $A$ is finite, $\#A$ is the number of elements contained in $A$. For any $z\in \Cx$ and $r>0$, we write $\Dx(z;r):=\{w\in \Cx: |w-z|<r\}$, $\overline{\Dx}(z;r):=\Cl(\Dx(z;r))$, $\Ex(z):=\Cx\setminus\overline{\Dx}(z;1)$, and $\overline{\Ex}(z):=\Cl(\Ex(z))$. We also write $\|z\|=\max\{\Re(z)|,|\Im(z)|\}$ for any $z\in\Cx$. More notation and terminology is defined throughout the text.

\section{Hurwitz Continued Fractions}
Let us define Hurwitz continued fraction. Denote by $\lfloor\cdot\rfloor:\RE\to\Za$ the usual floor function, that is, {$\lfloor x\rfloor$ is the largest integer less than or equal to $x$ for any real number $x$}. Define the function $[\cdot]:\Cx\to \Za[i]$ by
\[
\forall z\in\Cx \quad
[z]= \left\lfloor \Re(z) + \frac{1}{2}\right\rfloor + i\left\lfloor \Im(z) + \frac{1}{2}\right\rfloor,
\]
and the sets
\[
\mfF:=\{z\in\Cx:[z]=0\}, \quad \mfF^*:=\mfF\setminus \{0\}, \quad \mfF':=\mfF\setminus\QU(i).
\]
Let $a_1:\mfF^*\to \Za[i]$ be given by $a_1(z)= [z^{-1}]$ and $T:\mfF\to\mfF$ by
\[
\forall z\in\mfF \quad 
T(z) = 
\begin{cases}
\frac{1}{z} - \left[ \frac{1}{z}\right], \text{ if } z\neq 0, \\
0, \text{ if } z= 0.
\end{cases}
\]
As usual, $T^0:\mfF\to\mfF$ is the identity map and $T^{n}=T^{n-1}\circ T$ for $n\in\Na$. The function $T$ is a complex analogue of the usual Gauss map. 

Let $z$ be any complex number. Define the possibly finite sequence of Gaussian integers $(a_n(z))_{n\geq 0}$ as follows. First, $a_0(z):=[z]$ and, as long as $T^{n-1}(z - a_0(z)) \neq 0$, put 
\[
a_n(z) = a_1\left( T^{n-1} (z - a_0(z) )\right).
\]
The \textbf{Hurwitz continued fraction} (HCF) of $z$ is 
\[
[a_0;a_1,a_2,\ldots]:=  a_0 + \cfrac{1}{a_1 + \cfrac{1}{a_2 + \cfrac{1}{\ddots}}}.
\]
We shall refer to the sequence $(a_n(z))_{n\geq 0}$ as the \textbf{Hurwitz elements} of $z$ and we will omit the dependence on $z$ if there is no risk of ambiguity. 

Several properties of regular continued fractions hold too in the complex case, although their known proofs might be much more complicated (for example, the third point of proposition \ref{Propo01} below). Among the properties that remain true, we have that the HCF of any $z\in\QU(i)$ is finite and equals $z$. Also, when $z$ is an irrational complex number, that is $z\in\Cx':=\Cx\setminus \QU(i)$, the sequence $\sanu$ is infinite and 
\[
\lim_{n\to\infty} [a_0;a_1,\ldots,a_n] = z
\]
(see \cite{hurwitz} for a proof). 
\begin{propo01}\label{Propo01}
Let $z$ belong to $\Cx'$, let $\sanu$ be its Hurwitz elements and define the sequences of Gaussian integers $\sepn$, $\seqn$ by
\[
\begin{pmatrix}
p_{-1} & p_{-2} \\
q_{-1} & q_{-2}
\end{pmatrix}
=
\begin{pmatrix}
1 & 0 \\
0 & 1
\end{pmatrix}, \qquad
\forall n\in\Na \quad
\begin{pmatrix}
p_n\\
q_n
\end{pmatrix}
=
\begin{pmatrix}
p_{n-1} & p_{n-2} \\
q_{n-1} & q_{n-2}
\end{pmatrix}
\begin{pmatrix}
a_n\\
1
\end{pmatrix}.
\]
The following statements hold for all $n\in\Na_0$:
\begin{enumerate}[1.]
\item $q_{n}p_{n-1} - q_{n-1}p_{n}=(-1)^n$,
\item If $z_n = T^{n-1}(z-a_0)$, then {$z= \displaystyle\frac{p_{n-2}z_n+p_{n-1}}{q_{n-2}z_n+q_{n-1}}$},
\item $|q_n|< |q_{n+1}|$. \label{MonCont}
\item $\tfrac{p_n}{q_n}=[a_0;a_1,\ldots,a_n]$. \label{ApCont}
\end{enumerate}
\end{propo01}
The proofs of \ref{MonCont} and \ref{ApCont} can be found in \cite{hurwitz}. The rest of the theorem is proven in a broader context as Propositions 3.3 and 3.7 in \cite{daninog}. Proposition \ref{Propo01} also holds for finite HCF and any integer $n$ such that the corresponding objects are defined.

We shall now explore the geometry of $T$. Let $\iota:\Cx^*\to\Cx^*$ be the complex inversion, $\iota(z)=z^{-1}$, then
\begin{equation}\label{Eq-HCF-01}
\iota[\mfF] := \overline{\Ex}(-1)\cap \overline{\Ex}(i)\cap \Ex(1) \cap {\Ex(-i)},
\end{equation}
(see the paragraph on notation by then end of the introduction) and thus
\begin{equation}\label{Eq-HCF-02}
\forall z\in\Cx' \quad \forall n\in\Na \quad |a_n|\geq \sqrt{2}.
\end{equation}
Hence, if for any $\bfa=(a_1,\ldots,a_n)\in \Za[i]^n$ we define {the \textbf{cylinder}}
\[
\clC_n(\bfa):=\{z\in\mfF: a_1(z) = a_1,\ldots, a_n(z)=a_n\},
\]
we obtain from \eqref{Eq-HCF-02} that $\clC_1(a)=\vac$ for $a\in\{0,1,-1,i,-i\}$. Equation \eqref{Eq-HCF-01} along with some direct computations show that $T[\clC_1(a)]=\mfF$ when $|a|\geq \sqrt{8}$. In general, given $a\in\Za[i]$ with $|a|\in\{\sqrt{2},2,\sqrt{5}\}$, the set $T[\clC_1(a)]$ may assume fourteen different forms and, if we overlook the boundaries, it may assume only four different forms or any of their rotations by right angles (see Section 2 in \cite{eiitonana}). Moreover, if for any $n\in\Na$ and $\bfa\in\Za[i]^n$ we write 
\[
\mfF_n(\bfa):=T^n[\clC_n(\bfa)],
\]
an inductive argument shows that the set $\{\mfF_n(\bfa) : n\in\Na, \bfa\in \Za[i]^n\} \setminus\{\vac\}$ is finite and we may conclude that 
\[
\{ \leb\left(\mfF_n(\bfa)\right): n\in\Na, \bfa\in\Za[i]^n\}\setminus\{0\} 
\]
is non-empty and finite (cfr. Section 2 in \cite{eiitonana}).

It may happen that $\mfF_n(\bfa)\neq\vac$ but $\leb(\mfF_n(\bfa))=0$. Given $n\in\Na$, we say that $\bfa\in\Za[i]^n$ is \textbf{regular} if $\clC_n(\bfa)$ has non-empty interior and \textbf{irregular} otherwise. An infinite sequence in $\Za[i]$ is \textbf{regular} if all its prefixes are regular and irregular otherwise. Finally, a complex number is \textbf{regular} if the sequence of its HCF elements is regular and it is \textbf{irregular} in other case. Some properties follow directly from the definitions, inductive arguments, and our previous discussions. For example, if $\bfa\in\Za^n$ is irregular and $\clC_n(\bfa)\neq \vac$, then $\clC_n(\bfa)$ is contained {in some arc or in some line segment}. Therefore, the set of irregular complex numbers is null relative to the Lebesgue measure. Also, for any $n\in\Na$ and any regular $\bfa\in\Za[i]^n$ we have
\begin{equation}\label{Eq-HCF-03-05}
{\text{for some }} k\in\{0,1,2,3\} \quad 
i^k\inte \mfF_1({1-i})\subseteq \mfF_n(\bfa)\subseteq \mfF
\end{equation}
{(cfr. Lemma 1 in \cite{eiitonana})}. In what follows, we ignore irregular numbers, for they form a $\leb$-null set. However, for completeness sake, we exhibit an irregular sequence. If $n\in\Na_{\geq 2}$, then the set
\[
\iota\left[ \mfF_1(-1+i)\right] \cap \{z+ 1-in:z\in\mfF\} \neq\vac
\]
is a segment of the line determined by $\Re(z)=\frac{1}{2}$, so $(-1+i,1-in)$ is irregular.

We will not delve into the structure of the shift space associated to the dynamical system $(\mfF,T)$. We content ourselves with pointing out that it is rather complicated (see section 5.3 in \cite{hensley}). 

For any $z\in\Cx$, write {$\|z\|:=\max\{|\Re(z)|, |\Im(z)|\}$. Hence, for every $m\in\Na_{\geq 3}$ we have}
\begin{align*}
\#\left(\{ b\in\Za[i]: \|b\|=m\} \cap \iota\left[\mfF\right]\right) &= 8m, \nonumber\\
\#\left(\{ b\in\Za[i]: \|b\|=m\} \cap \iota\left[\inte\mfF_1({1-i})\right]\right) &= 2m-1. \nonumber
\end{align*}
{For $m=2$, the first set contains fourteen elements while the second set has three elements. In} view of \eqref{Eq-HCF-03-05}, we can conclude the next estimate:
\begin{lem01}\label{Le-HCF-01}
Let $n\in\Na$ and $\bfa\in\Za[i]^n$ be such that $\clC_n(\bfa)$ is regular, then
\begin{equation}\label{Eq-HCF-04}
\forall m\in\Na_{\geq 2}\qquad
\#\; \left\{ b\in\Za[i]: \|b\|=m\right\} \cap \iota\left[\mfF_n(\bfa)\right] \asymp m.
\end{equation}
\end{lem01}

\begin{lem01}\label{Le-HCF-bd}
Every $z=[0;a_1,a_2,\ldots]\in\mfF'$ satisfies {
\begin{equation}\label{Le01-Eq01}
\forall n\in\Na_0\qquad
\left| z - \frac{p_n}{q_n}\right|\leq  \frac{4+2\sqrt{2}}{|a_{n+1}||q_n|^2 }.
\end{equation}}
\end{lem01}
\begin{proof}
Take $z=[0;a_1,a_2,\ldots]\in\mfF'$ and $n\in\Na$. Write $z_j=[0;a_j,a_{j+1},a_{j+2},\ldots]$ for every $j\in\Na$ and let $\sepn$, $\seqn$ be as in Proposition \ref{Propo01}, then
\[
\left| z -\frac{p_n}{q_n}\right| = \left| \frac{p_n z_{n+1}^{-1} + p_{n-1}}{q_nz_{n+1}^{-1} + q_{n-1}} - \frac{p_n}{q_n}\right| = \frac{1}{|q_n|^2\left| z_{n+1}^{-1} + \frac{q_{n-1}}{q_n}\right|}.
\]
(see Proposition \ref{Propo01}). Thus, we can conclude \eqref{Le01-Eq01} if we obtain an absolute constant $\kappa_1>0$ such that
\begin{equation}\label{Le01-Eq02}
\forall n\in\Na \quad
\left| z_{n+1}^{-1} + \frac{q_{n-1}}{q_n} \right| \geq \kappa_1 |a_{n+1}|.
\end{equation}
Take $n\in\Na$. Assume that $|a_{n+1}|\geq 2$. Put
\[
{\kappa_1}= 1- \frac{1}{2}\left(1+ \frac{1}{\sqrt{2}}\right) = \frac{2-\sqrt{2}}{4}.
\]
Since $|q_{n-1}|<|q_n|$ and $|z_{n+2}|\leq 1/\sqrt{2}$, we have
\[
\left| {z_{n+1}^{-1}} + \frac{q_{n-1} }{q_n} \right|
=
\left| a_{n+1} + {z_{n+2}} + \frac{q_{n-1}}{q_n} \right|
\geq 
|a_{n+1}| - \left( 1 + \frac{1}{\sqrt{2}} \right) 
\geq 
{\kappa_1}|a_{n+1}|.
\]
Now, assume that $|a_{n+1}|=\sqrt{2}$. Since $\Cl(\mfF) \subseteq \Dx(0;1)$, the closed set $\iota[\Cl (\mfF)]$ is entirely contained in the complement of the compact set $\overline{\Dx}(0;1)$, so
\[
\kappa_2:=\inf\left\{ |w-w'|: w\in\iota[\Cl (\mfF)], \; w'\in\overline{\Dx}(0;1)\right\}>0,
\]
{and hence}
\[
\left| {z_{n+1}^{-1}} + \frac{q_{n-1}}{q_n}\right| \geq \kappa_2 = \frac{\kappa_2}{\sqrt{2}}|a_{n+1}|.
\]
{It is not hard to show that $\tfrac{\kappa_2}{\sqrt{2}} >\kappa_1$, so \eqref{Le01-Eq02} also holds in this case. Because $\kappa_1^{-1}=4+2\sqrt{2}$, the lemma is proven.}
\end{proof}
Once again, let us introduce some notation. For any $n\in\Na$, $\bfa=(a_1,\ldots,a_n)\in\Za[i]^n$, $m\in\Na$, and $\bfb=(b_1,\ldots,b_m)\in\Za[i]^m$, write $\bfa\bfb:=(a_1,\ldots,a_n,b_1,\ldots,b_m)$, and define the sets $\Omega_{R}(n):=\{\bfc\in\Za[i]^n: \bfc \text{ is regular} \}$ and
\[
R(\bfa;m):=\left\{ \bfc\in\Za[i]^m: \bfa\bfc\in \Omega_R(n+m)\right\}.
\]
\begin{lem01}\label{Le-HCF-02}
If $n\in\Na$, $\bfa\in \Omega_R(n)$, and $b\in R(\bfa;1)$, then
\[
\frac{\leb\left( \clC_{n+1}(\bfa b) \right)}{\leb\left( \clC_{n}(\bfa) \right)} \asymp \frac{1}{|b|^4}.
\]
\end{lem01}
{Lemma \ref{Le-HCF-02} is the complex version of a well known estimate for regular continued fractions (see Equation (57) in \cite{khin} for a precise statement and its proof).}
\begin{proof}
Let $n,\bfa,b$ be as in the statement and take $\sepn$ and $\seqn$ as in Proposition \ref{Propo01}. The restriction of $T^n$ to $\clC_n(\bfa)$ is a bijection onto $\mfF_n(\bfa)$ with inverse $t_{\bfa,n}:\mfF_n(\bfa)\to\clC_n(\bfa)$ 
\[
\forall z\in\mfF_n(\bfa) \quad
t_{\bfa,n}(z) = \frac{p_{n-1}z + p_n}{q_{n-1}z + q_n} = \frac{p_{n-1}}{q_{n-1}} + \frac{{(-1)^{n-1}}}{q_{n-1}(q_{n-1} z + q_n)}.
\]
Since $|q_{n-1}|< |q_n|$ and $\mfF_n(\bfa)\subseteq \Dx\left(0; 1/\sqrt{2}\right)$, $t_{\bfa,n}$ is holomorphic and
\[
\forall z\in\mfF_n(\bfa) \quad t_{\bfa,n}'(z) = \frac{{(-1)^{n}} }{(q_{n-1}z+ q_n)^2}.
\]
When we identify $\Cx$ with $\RE^2$, the function $t_{\bfa,n}$ is differentiable on $\mfF_n(\bfa)$ and, if $Dt_{\bfa,n}(x,y)$ denotes its real derivative at the point $(x,y)$, we get {for all $z=x+iy$ in $\mfF_n(\bfa)$}
\[
|\det D t_{\bfa,n}(x,y)| = {|t_{\bfa,n}'(z)| }=\frac{1}{|q_{n-1}z + q_n|^4} = \frac{1}{|q_n|^4 \left| \frac{q_{n-1}}{q_n}z + 1 \right|^4}
\]
(see p. 33 in \cite{lang-ca}) and, since {every $z\in\mfF$ satisfies}
\[
1-\frac{1}{\sqrt{2}} \leq 1-|z|\leq 1- \left| \frac{q_{n-1}}{q_n} \right| |z|
\leq 
\left| \frac{q_{n-1}}{q_n}z +1 \right| 
\leq \left|\frac{q_{n-1}}{q_n}z\right| + 1 < 2,
\]
we may conclude that $|\det D t_{\bfa,n}|\asymp |q_n|^{-4}$. Therefore, the Theorem of Change of Variable (Theorem 3.7.1 in \cite{bogachev}, Vol. 1., p. 194) {implies} that every Borel set $B\subseteq \Cx$ satisfies
\[
\leb\left(t_{\bfa,n}[B\cap \mfF_n(\bfa)]\right)\asymp \frac{\leb\left(B\cap \mfF_n(\bfa)\right)}{|q_n|^4}.
\]
In particular, choosing $B=\mfF$ {and using \eqref{Eq-HCF-03-05}}, we get $\leb(\clC_{n}(\bfa))\asymp |q_n|^{-4}$. {If we apply a similar argument on $\clC_{n+1}(\bfa b)$}, we arrive at
\[
\frac{\leb\left( \clC_{n+1}(\bfa b) \right)}{\leb\left( \clC_{n}(\bfa) \right)}
\asymp
\frac{|q_n |^4}{|bq_n +q_{n-1}|^4} = \frac{1}{\left| b+ \frac{q_{n-1}}{q_n}\right|^4}.
\]
Finally, by $|q_{n-1}|<|q_n|$ and $\sqrt{2}\leq |b|$ (see \eqref{Eq-HCF-02}), the constant ${\kappa_3}=4(\sqrt{2}-1)^{-4}$ verifies
\[
\frac{1}{2^4|b|^4} \leq \frac{1}{(|b|+1)^4}
\ll 
\frac{\leb\left( \clC_{n+1}(\bfa b) \right)}{\leb\left( \clC_{n}(\bfa) \right)}
\ll
\frac{1}{(|b|-1)^4}
\leq
\frac{{\kappa_3}}{|b|^4}. \qedhere
\]
\end{proof}

\section{Proof of the Main Theorem}

Let $\bfu=(u_n)_{n\geq 1}$ be a sequence in $\RE_{>0}$. Instead of working with $E(\bfu)$, we will use more manageable sets. For any $\kappa>0$ define
\[
E^{\infty}(\kappa \bfu):= \left\{ z=[0;a_1,a_2,\ldots] \in\mfF: \kappa u_n \leq \|a_n\| \quad\text{for infinitely many } n\in\Na \right\}
\]
and $E^{\infty} (\bfu):=E^{\infty} (1\bfu)$. Note that
\begin{equation}\label{Eq-Cont}
E^{\infty} (\bfu) \subseteq E(\bfu)\subseteq E^{\infty}\left(\frac{1}{\sqrt{2}} \bfu\right),
\end{equation}
because $\|z\| \leq |z| \leq \sqrt{2}\|z\|$ for any $z\in\Cx$.
\begin{lem01}\label{Le01}
If $\sum_n u_n^{-2}<+\infty$, then $\leb(E^{\infty}(\kappa\bfu))=0$ for every $\kappa>0$.
\end{lem01}
\begin{proof}
Define the sequence of sets $(E_n^{\infty})_{n\geq 1}$ by
\begin{equation}\label{Eq-Def-En}
\forall n\in\Na \qquad
E_n^{\infty}:=\left\{ [0;a_1,a_2,\ldots] \in\mfF': u_n\leq \|a_n\| \right\}
\end{equation}
For any $k\in\Na$ and $\bfa\in\Omega_R(k)$, {Lemmas \ref{Le-HCF-01} and \ref{Le-HCF-02} give}
\begin{align*}
\leb\left(E_{k+1}^{\infty}\cap \clC_k(\bfa)\right) &= \sum_{\|b\| \geq u_{k+1}} \leb\left( \clC_{k+1}(\bfa b)\right) \nonumber\\
&\asymp \sum_{\|b\| \geq u_{k+1}} \frac{\leb\left( \clC_{k}(\bfa)\right)}{|b|^4} \nonumber\\
&\ll \leb\left( \clC_k(\bfa) \right) \sum_{n\geq u_{k+1}} n^{-3} {\leq \frac{\leb\left( \clC_k(\bfa)\right)}{(u_{k+1}-1)^2},}
\end{align*}
{Note that there is no loss of generality in assuming that $u_{k+1}>1$.} Letting $\bfa$ run along $\Za[i]^k$, we arrive at $\leb(E_{k+1}^{\infty}){\ll (u_{k+1}-1)^{-2}}$. Hence, by the Borel-Cantelli lemma, $\leb\left(E^{\infty}(\bfu)\right)=0$. Finally, for any $\kappa>0$, the convergence of $\sum_k u_k^{-2}$ implies that of $\sum_k (\kappa u_k)^{-2}$ and the previous argument shows $\leb\left( E^{\infty}(\kappa\bfu)\right)=0$.
\end{proof}
\begin{lem01}\label{Le02}
If $\sum_{n\geq 1} u_n^{-2}=+\infty$, then $\leb\left( E^{\infty}(\kappa\bfu)\right)=1$ for every $\kappa>0$.
\end{lem01}
\begin{proof}
{As in the proof of Lemma \ref{Le01},} it is enough to show the result for $\kappa=1$. Let $(E_n^{\infty})_{n\geq 1}$ be as in \eqref{Eq-Def-En} and define
\[
\forall n\in\Na \qquad
B_n:= \bigcap_{j\geq n} \mfF'\setminus E_j^{\infty}.
\]
Note that $(B_n)_{n\geq 1}$ is increasing and
\[
\mfF'\setminus E^{\infty}(\bfu) = \bigcup_{n\in\Na} B_n.
\]
Take $k\in\Na$, $\bfa\in \Omega_R(k)$ and let $0<c_1<1$ be a constant (independent of $k$ and $\bfa$) such that for any $k\in\Na$ and any $\bfa\in \Omega_R(k)$ we have
\[
\forall b\in R(\bfa;1) \qquad \leb\left(\clC_{k+1}(\bfa b)\right) \geq \frac{c_1}{|b|^4} \leb\left({\clC_{k}}(\bfa)\right)
\]
(cfr. Lemma \ref{Le-HCF-02}). Let $0<c_2<1$ be such that for every $M>0$
\[
\sum_{j\geq M} \frac{1}{j^3} > c_2 \frac{1}{(M+1)^2}.
\]
Therefore, {after choosing a suitable constant $0<c_3<1$ with the aid of Lemma \ref{Le-HCF-01} and writing $c=c_1c_2c_3$}, we have
\[
\sum_{\substack{\|b\| \geq M \\ b\in R(\bfa;1)}} \leb\left({\clC_{k+1}}(\bfa b)\right)> c_1{c_3} \leb\left( \clC_k(\bfa)\right) \sum_{j\geq M} \frac{1}{j^3} > c \frac{\leb\left(\clC_k(\bfa)\right)}{(M+1)^2},
\]
and thus
\begin{equation}\label{Eq-LowBow}
\sum_{\|b\| \leq M} \leb\left( \clC_{k+1}(\bfa b)\right) \leq \left( 1- \frac{c}{(M+1)^2} \right) \leb\left( \clC_{k}(\bfa)\right)
\end{equation}
Define for each $n\in\Na$ the set 
\[{
F(\bfa,k,n)=\bigcup_{\bfb} \clC_{k+n}(\bfa\bfb), }
\]
{where $\bfb=(b_1,\ldots,b_n)$ runs along the $n$-tuples that belong to $R(\bfa;n)$ and that satisfy $\|b_j\|<u_{k+j}$ for all $j\in\{1,\ldots,n\}$.} Then, applying \eqref{Eq-LowBow} recursively,
\begin{equation}\label{Eq-MT-03}
\forall n\in\Na_{\geq 2} \quad
\leb\left( F(\bfa,k,n)\right) \leq \leb\left( F(\bfa,k,1)\right) \prod_{j=2}^n \left( 1- \frac{c}{(1+u_{k+j})^2}\right).
\end{equation}
Because $1+x\leq e^x$ holds for $x\in\RE$, the equality
\[
B_{k+1}\cap \clC_k(\bfa) = \bigcap_{n\in\Na} F(\bfa,k,n)
\]
and \eqref{Eq-MT-03} imply $\leb\left( B_{k+1}\cap \clC_k(\bfa)\right)=0$. Taking the union along $\bfa\in \Za[i]^k$, we obtain $\leb(B_{k+1})=0$ and hence
\[
\leb\left(\mfF'\setminus E^{\infty}(\bfu)\right)= \leb\left( \bigcup_{k\in\Na} B_k\right) =0. \quad \qedhere
\]
\end{proof}
\begin{proof}[Proof of Theorem \ref{TEOBBC}]
On the one hand, if $\sum_{n}u_n^{-2}<+\infty$, then Lemma \ref{Le01} and \eqref{Eq-Cont} give
\[
\leb\left(E(\bfu)\right)\leq  \leb\left(E^{\infty}\left(\frac{1}{\sqrt{2}}\bfu\right)\right)=0.
\]
On the other hand, $\sum_n u_n^{-2}=+\infty$ and Lemma \ref{Le02} imply $\leb(\mfF\setminus E^{\infty}(\bfu))=0$ and
\[
1 
= \leb\left(E^{\infty} (\bfu) \right) 
\leq \leb\left( E(\bfu) \right) 
\leq \leb(\mfF)=1. \qedhere
\] 
\end{proof}

\subsection{Proof of Corollary \ref{COROKHIN}}
The well known existence of a Borel probability measure on $\mfF$ which is equivalent to $\leb$ and $T$-ergodic yields the following result (see propositions 1 and 2 in \cite{nakada87}):
\begin{propo01}\label{PROP-LK-CONSTANT}
There exists a positive constant $B>0$ such that for $\leb$-almost every $z\in\mfF$ the sequence $\seqn$ defined as in Proposition \ref{Propo01} satisfies
\[
\lim_{n\to\infty} \frac{1}{n} \log|q_n| = B.
\]
\end{propo01}

Take any function $\psi:\RE_{>0}\to\RE_{>0}$. As in the regular continued fraction context, we only need $x\mapsto x^2\psi(x)$ to be non-increasing in the divergence case. Observe that, because $z-[z]$ belongs to $\mfF$ for any $z\in\Cx$, the set $K(\psi)$ is the union of Gaussian integral translates of $K(\psi)\cap\mfF$. Therefore, it suffices to show Corollary \ref{COROKHIN} for $K(\psi)\cap\mfF$.

\begin{proof}
Assume that $\sum_n n^{3}\psi(n)^2<+\infty$. Define for each $n\in\Na$ the sets
\begin{align*}
Z_n =\{z\in\Za[i]:\|z\|=n\}, &\quad C_n=\{z\in\Za[i]: \|z\| < n\}, \nonumber\\
K_n = &\bigcup_{\substack{p\in C_n\\q\in Z_n }} \overline{\Dx}\left(\frac{p}{q},\psi(|q|) \right)\cap \mfF. \nonumber
\end{align*}
Then, we have that
\[
\leb(K_n)  \ll \#Z_n\, \#C_n\psi(n)^2 \ll n^3\psi(n)^2,
\]
so $\sum_n \leb(K_n)<+\infty$ and, by the Borel-Contelli Lemma, $\leb(K(\psi)\cap \mfF)=0$. Then, the discussion prior to the proof implies $\leb(K(\psi))=0$.

Now assume that $\sum_n n^{3}\psi(n)^2=+\infty$ and that $x\mapsto x^2\psi(x)$ is non-increasing for $x>0$. Then, $y\mapsto y^3\psi(y)^2$ is non-increasing for $y>1$. Let $D$ be a positive number strictly larger than $B$ (given as in proposition \ref{PROP-LK-CONSTANT}) and define $\tilde{\Phi}:\RE_{>0}\to \RE_{>0}$ by
\[
\forall x\in \RE_{>0} \qquad \tilde{\Phi}(x) =e^{2Dx} \psi(e^{Dx}).
\]
Hence, $\tilde{\Phi}$ is non-increasing and for $0<t<T$ we have
\[
\int_t^T \tilde{\Phi}^2(x)\md x 
=
\int_t^T e^{4Dx} \psi(e^{Dx})^2 \md x 
=
\frac{1}{D}\int_{e^{Dt}}^{e^{DT}} y^3 \psi(y)^2\md y
\]
and we conclude $\sum_{n}\tilde{\Phi}(n)^2=+\infty$ by letting $T\to + \infty$.

{Write $\kappa=4+2\sqrt{2}$} and define $\Phi=\kappa^{-1}\tilde{\Phi}$. By the main theorem, for almost all $z=[0;a_1,a_2,\ldots]\in\mfF$ the inequality $|a_{n+1}|\geq \Phi(n)^{-1}$ holds infinitely often. Also, because of Proposition \ref{PROP-LK-CONSTANT}, we have $|q_n|<e^{Dn}$ for sufficiently large $n$ almost everywhere. Therefore, {by Lemma \ref{Le-HCF-bd}}, for almost every $z\in\mfF$ there are infinitely many $n\in\Na$ satisfying
\[
\left| z - \frac{p_n}{q_n}\right| \leq \frac{\kappa}{|a_{n+1}||q_n|^2} \leq \frac{\tilde{\Phi}(n)}{|q_n|^2} \leq \frac{\tilde{\Phi}\left(D^{-1}\log|q_n|\right)}{|q_n|^2} =  \psi(|q_n|). \qedhere
\]
\end{proof}
\section{Acknowledgements}
{I thank Yann Bugeaud, Mauricio González Soto, Mumtaz Hussain, Simon Kristensen, and the referees for their helpful comments. I especially thank M. Hussain for suggesting this problem.}

\bibliographystyle{acm}
\bibliography{referencias}

\begin{thebibliography}{10}

\bibitem{bdv}
{\sc Beresnevich, V., Dickinson, D., and Velani, S.}
\newblock Measure theoretic laws for lim sup sets.
\newblock {\em Mem. Amer. Math. Soc. 179}, 846 (2006), x+91.

\bibitem{bogachev}
{\sc Bogachev, V.~I.}
\newblock {\em Measure theory. {V}ol. {I}, {II}}.
\newblock Springer-Verlag, Berlin, 2007.

\bibitem{bugeaud-livre}
{\sc Bugeaud, Y.}
\newblock {\em Approximation by algebraic numbers}, vol.~160 of {\em Cambridge
  Tracts in Mathematics}.
\newblock Cambridge University Press, Cambridge, 2004.

\bibitem{daninog}
{\sc Dani, S.~G., and Nogueira, A.}
\newblock Continued fractions for complex numbers and values of binary
  quadratic forms.
\newblock {\em Trans. Amer. Math. Soc. 366}, 7 (2014), 3553--3583.

\bibitem{dokr}
{\sc Dodson, M.~M., and Kristensen, S.}
\newblock Hausdorff dimension and {D}iophantine approximation.
\newblock In {\em Fractal geometry and applications: a jubilee of {B}eno\^{\i}t
  {M}andelbrot. {P}art 1}, vol.~72 of {\em Proc. Sympos. Pure Math.} Amer.
  Math. Soc., Providence, RI, 2004, pp.~305--347.

\bibitem{eiitonana}
{\sc Ei, H., Ito, S., Nakada, H., and Natsui, R.}
\newblock On the construction of the natural extension of the {H}urwitz complex
  continued fraction map.
\newblock {\em Monatsh. Math. 188}, 1 (2019), 37--86.

\bibitem{HeXi}
{\sc He, Y., and Xiong, Y.}
\newblock The difference between the {H}urwitz continued fraction expansions of
  a complex number and its rational approximations.
\newblock {\em Fractals\/} (2021), to appear.

\bibitem{hensley}
{\sc Hensley, D.}
\newblock {\em Continued fractions}.
\newblock World Scientific Publishing Co. Pte. Ltd., Hackensack, NJ, 2006.

\bibitem{hurwitz}
{\sc Hurwitz, A.}
\newblock \"{U}ber die {E}ntwicklung complexer {G}r\"{o}ssen in
  {K}ettenbr\"{u}che.
\newblock {\em Acta Math. 11}, 1-4 (1887), 187--200.

\bibitem{hussain17}
{\sc Hussain, M.}
\newblock A {K}hintchine-{G}roshev type theorem in absolute value over complex
  numbers.
\newblock {\em New Zealand J. Math. 47\/} (2017), 57--67.

\bibitem{khin}
{\sc Khinchin, A.~Y.}
\newblock {\em Continued fractions}, {R}ussian~ed.
\newblock Dover Publications, Inc., Mineola, NY, 1997.
\newblock With a preface by B. V. Gnedenko, Reprint of the 1964 translation.

\bibitem{lang-ca}
{\sc Lang, S.}
\newblock {\em Complex analysis}, fourth~ed., vol.~103 of {\em Graduate Texts
  in Mathematics}.
\newblock Springer-Verlag, New York, 1999.

\bibitem{leveque}
{\sc LeVeque, W.~J.}
\newblock Continued fractions and approximations in {$k(i)$}. {I}, {II}.
\newblock {\em Nederl. Akad. Wetensch. Proc. Ser. A. {\bf 55} = Indagationes
  Math. 14\/} (1952), 526--535, 536--545.

\bibitem{nakada87}
{\sc Nakada, H.}
\newblock On metric {D}iophantine approximation of complex numbers, complex
  continued fractions.
\newblock In {\em S\'{e}minaire de {T}h\'{e}orie des {N}ombres, 1987--1988
  ({T}alence, 1987--1988)}. Univ. Bordeaux I, Talence, [1988?], pp.~Exp. No.
  45, 10.

\bibitem{nogueira}
{\sc Nogueira, A.}
\newblock The {B}orel-{B}ernstein theorem for multidimensional continued
  fractions.
\newblock {\em J. Anal. Math. 85\/} (2001), 1--41.

\bibitem{sullivan}
{\sc Sullivan, D.}
\newblock Disjoint spheres, approximation by imaginary quadratic numbers, and
  the logarithm law for geodesics.
\newblock {\em Acta Math. 149}, 3-4 (1982), 215--237.

\end{thebibliography}

\end{document}